\documentclass[12pt]{article}
 \usepackage{amssymb,amsmath,amsthm,epsfig,graphics,latexsym}
 
\begin{document}

 \begin{center}
 {\bf On Relatively Prime Subsets and Supersets}
 \vskip 20pt
 {\bf Mohamed El Bachraoui\footnote{Supported by RA at UAEU, grant: 02-01-2-11/09}} \\
  \emph{Dept. Math. Sci., United Arab Emirates University, P.O.Box 17551, Al-Ain, UAE} \\
 {\tt melbachraoui@uaeu.ac.ae}\\

 \vskip 10pt
 \end{center}
 \vskip 30pt
 \vskip 30pt

 \centerline{\bf Abstract}

 \noindent
 A nonempty finite set of positive integers $A$ is relatively prime if $\gcd(A)=1$ and it is
 relatively prime to $n$ if $\gcd(A\cup \{n\})=1$.
 The number of nonempty subsets of $A$ which are relatively prime to $n$ is $\Phi(A,n)$
 and the number of such subsets of cardinality $k$ is $\Phi_k(A,n)$. Given positive integers
 $l_1$, $l_2$, $m_2$, and $n$ such that $l_1\leq l_2\leq m_2$ we give
 $\Phi( [1,m_1]\cup [l_2, m_2],n)$ along with $\Phi_k( [1,m_1]\cup [l_2, m_2],n)$.
 Given positive integers $l, m$, and $n$ such that $l\leq m$
 we count for any subset $A$ of $\{l,l+1,\ldots,m\}$ the number of its supersets in $[l,m]$ which are
 relatively prime and we count the number of such supersets which are relatively prime to $n$.
 Formulas are also obtained for corresponding supersets having fixed cardinalities.
 Intermediate consequences include a formula for the number of relatively prime sets
 with a nonempty intersection with some fixed set of positive integers.

 \pagestyle{myheadings}

 \thispagestyle{empty}
 \baselineskip=15pt
 \vskip 30pt
 \noindent
  \\ \\
  \noindent
  \textbf{Subject Class:}\quad 11A25, 11B05, 11B75.

 \section*{\normalsize 1. Introduction}

 Throughout let $k, l, m , n$ be positive integers such that $l \leq m$, let $[l,m] = \{l,l+1,\ldots,m\}$,
 let $\mu$ be the M\"obius function, and let $\lfloor x \rfloor$
  be the floor of $x$. If $A$ is a set of integers and $d\not= 0$, then
  $\frac{A}{d}= \{ a/d:\ a \in A\}$.
 A nonempty set of positive integers $A$ is called \emph{relatively prime} if $\gcd(A)=1$ and it is
 called \emph{relatively prime to $n$} if $\gcd(A\cup \{n\}) = \gcd(A,n) = 1$. Unless otherwise specified
 $A$ and $B$ will denote nonempty sets of positive integers.
 We will need the following basic identity on binomial coefficients stating that for nonnegative integers
 $L\leq M \leq N$
 \begin{equation}\label{binomial}
  \sum_{j=M}^{N}\binom{j}{L} = \binom{N+1}{L+1}-\binom{M}{L+1}.
 \end{equation}


 \noindent
 {\bf Definition 1.}
 Let
  \[
  \begin{split}
  \Phi(A,n) &= \# \{X\subseteq A:\ X\not=
  \emptyset\ \text{and\ } \gcd(X,n) = 1 \}, \\
   \Phi_k (A,n) &= \# \{X\subseteq A:\ \# X=
  k \ \text{and\ } \gcd(X,n) = 1 \}, \\
  f(A) &= \# \{X\subseteq A:\ X\not=
  \emptyset\ \text{and\ } \gcd(X) = 1 \}, \\
   f_k (A) &= \# \{X\subseteq A:\ \# X=
  k \ \text{and\ } \gcd(X) = 1 \}.
  \end{split}
  \]

 \noindent
 Nathanson in \cite{Nathanson} introduced $f(n)$, $f_k(n)$, $\Phi(n)$, and $\Phi_k(n)$
 (in our terminology $f([1,n])$, $f_k([1,n])$, $\Phi([1,n],n)$,
 and $\Phi_k([1,n],n)$ respectively) and gave their formulas along with asymptotic estimates.
 Formulas for $f([m,n])$, $f_k([m,n])$, $\Phi([m,n],n)$, and $\Phi_k([m,n],n)$ are found in \cite{ElBachraoui1, Nathanson-Orosz}
  and formulas for $\Phi([1,m],n)$ and $\Phi_k([1,m],n)$ for $m\leq n$ are obtained in
  \cite{ElBachraoui2}.
 %
 Recently Ayad and Kihel in \cite{Ayad-Kihel2} considered phi functions for sets which are in arithmetic progression
 and obtained the following more general formulas for $\Phi([l,m],n)$ and $\Phi_k ([l,m],n)$.

 \noindent
 {\bf Theorem 1.}\ 
 We have
 \[
 \begin{split}
  \text{(a)\quad } &\  \Phi([l,m],n) = \sum_{d|n}\mu(d) 2^{\lfloor m/d \rfloor- \lfloor (l-1)/d \rfloor}, \\
  \text{(b)\quad } &\ \Phi_k ([l,m],n) = \sum_{d|n} \mu(d)
   \binom{\lfloor m/d \rfloor- \lfloor (l-1)/d \rfloor}{k}.
  \end{split}
 \]
 \section*{\normalsize 2. Relatively prime subsets for $[1,m_1]\cup [l_2,m_2]$}
 If $[1,m_1]\cap[l_2,m_2]= \emptyset$, then phi functions for $[1,m_1]\cup[l_2,m_2]= [1,m_2]$ are
  obtained by Theorem 1. So we may assume that $1 \leq m_1 < l_2 \leq m_2$.
 %

 \noindent
 {\bf Lemma 1.}\label{lem:psi}
 Let
 \[
 \Psi(m_1,l_2,m_2, n)= \# \{X \subseteq [1,m_1]\cup [l_2,m_2]:\ l_2\in X\
 \text{and\ } \gcd(X,n)=1 \},
 \]
 \[
 \Psi_k(m_1,l_2,m_2,n)=\# \{X \subseteq [1,m_1]\cup [l_2,m_2]:\ l_2\in X,\ |X| = k,\
 \text{and\ } \gcd(X,n)=1 \}.
 \]
 Then
 \[
 \text{(a)\quad } \Psi(m_1,l_2,m_2, n) = \sum_{d|(l_2,n)}\mu(d)
  2^{\lfloor m_1/d \rfloor + \lfloor m_2/d \rfloor- l_2/d}, \]
\[
 \text{(b)\quad } \Psi_k (m_1,l_2,m_2, n) = \sum_{d|(l_2,n)}\mu(d)
 \binom{\lfloor m_1/d \rfloor + \lfloor m_2/d \rfloor - l_2/d }{k-1}.
 \]

\begin{proof}
 (a) Assume first that $m_2\leq n$. Let $\mathcal{P}(m_1,l_2,m_2)$ denote the set of subsets of
 $[1,m_1]\cup[l_2,m_2]$
 containing $l_2$ and let $\mathcal{P}(m_1,l_2,m_2,d)$ be the set of subsets $X$ of
 $[1,m_1]\cup[l_2,m_2]$ such that
 $l_2\in X$ and $\gcd(X,n) = d$. It is clear that the set $\mathcal{P}(m_1,l_2,m_2)$
 of cardinality $2^{m_1+m_2-l_2}$  can be
 partitioned using the equivalence relation of having the same $\gcd$ (dividing $l_2$ and $n$).
 Moreover, the mapping
 $A \mapsto \frac{1}{d} X$
 is a one-to-one correspondence between
 $\mathcal{P}(m_1,l_2,m_2,d)$ and the
 set of subsets $Y$ of
 $[1, \lfloor m_1/d \rfloor ]\cup [l_2/d,\lfloor m_2/d \rfloor]$
 such that $l_2/d \in Y$ and $\gcd(Y,n/d)= 1$. Then
 \[
  \# \mathcal{P}(m_1,l_2,m_2,d) =
  \Psi(\lfloor m_1/d \rfloor,l_2 /d,\lfloor m_2/d \rfloor,n/d).
 \]
 Thus
 \[
 2^{m_1+m_2-l_2} = \sum_{d|(l_2,n)} \# \mathcal{P}(m_1,l_2,m_2,d)=
 \sum_{d|(l_2,n)} \Psi (\lfloor m_1/d \rfloor,l_2 /d,\lfloor m_2/d \rfloor,n/d),
 \]
 which by the M\"obius inversion formula extended to multivariable functions
  \cite[Theorem 2]{ElBachraoui1} is equivalent to
 \[
 \Psi(m_1,l_2,m_2,n) = \sum_{d|(l_2,n)}\mu(d)
 2^{\lfloor m_1/d \rfloor + \lfloor m_2/d \rfloor - l_2/d}.
 \]
 Assume now that $m_2 >n$ and let $a$ be a positive integer such that $m_2 \leq n^a$.
 As $\gcd(X,n)=1$ if and only if $\gcd(X,n^a)=1$ and $\mu(d) =0$ whenever $d$ has a nontrivial
 square factor, we have
 \[
 \begin{split}
 \Psi(m_1,l_2,m_2,n) &= \Psi(m_1,l_2,m_2,n^a) \\
 &=
 \sum_{d|(l_2,n^a)}\mu(d)
 2^{\lfloor m_1/d \rfloor + \lfloor m_2/d \rfloor - l_2/d} \\
 &=
 \sum_{d|(l_2,n)}\mu(d)
 2^{\lfloor m_1/d \rfloor + \lfloor m_2/d \rfloor - l_2/d}.
 \end{split}
 \]

 (b) For the same reason as before, we may assume that $m_2 \leq n$.
  Noting that the correspondence
 $X\mapsto \frac{1}{d} X$ defined above preserves the cardinality and using an argument similar to
 the one in part (a), we obtain the following identity
 \[
  \binom{m_1+m_2-l_2}{k-1}=
  \sum_{d|(l_2,n)}\Psi_k (\lfloor m_1/d \rfloor,l_2 /d,\lfloor m_2/d \rfloor, n/d )
  \]
  which by the M\"obius inversion formula \cite[Theorem 2]{ElBachraoui1} is equivalent to
  \[
   \Psi_k (m_1,l_2,m_2,n) =
   \sum_{d|(l_2,n)}\mu(d)\binom{\lfloor m_1/d \rfloor + \lfloor m_2/d \rfloor - l_2/d }{k-1},
 \]
 as desired.
\end{proof}

 \noindent
 {\bf Theorem 2.}\label{thm:main2}
 We have
  \[
 \begin{split} \text{(a)\quad }
  \Phi([1,m_1]\cup [l_2,m_2],n) &= \sum_{d|n}\mu(d)
   2^{\lfloor \frac{m_1}{d} \rfloor +\lfloor \frac{m_2}{d} \rfloor - \lfloor\frac{l_2 -1}{d} \rfloor}, \\
  \text{(b)\quad }
  \Phi_k ([1,m_1]\cup [l_2,m_2],n) &=
   \sum_{d|n} \mu(d)
   \ \binom{\lfloor \frac{m_1}{d} \rfloor +\lfloor \frac{m_2}{d} \rfloor - \lfloor\frac{l_2 -1}{d} \rfloor}{k}.
   \end{split}
  \]
%
 \begin{proof}
 (a) Clearly
 \begin{equation}\label{help1}
 \begin{split}
 \Phi([1,m_1]\cup [l_2,m_2],n)
 & =
 \Phi([1,m_1]\cup [l_2 -1,m_2],n) - \Psi(m_1,l_2 -1,m_2,n) \\
 &=
 \Phi([1,m_1]\cup [m_1+1,m_2],n) - \sum_{i=m_1 +1}^{l_2 -1}\Psi(m_1,i,m_2,n) \\
 &=
 \Phi([1,m_2] - \sum_{i=m_1 +1}^{l_2 -1}\Psi(m_1,i,m_2,n) \\
 &=
 \sum_{d|n} \mu(d) 2^{\lfloor m_2/d \rfloor} - \sum_{i=m_1 +1}^{l_2 -1}\sum_{d|(n, i)}
 \mu(d) 2^{\lfloor \frac{m_1}{d} \rfloor +\lfloor \frac{m_2}{d} \rfloor - \frac{i}{d}},
 \end{split}
 \end{equation}
 where the last identity follows by Theorem 1 for $l=1$ and Lemma 1.
 Rearranging the last summation in (\ref{help1}) gives
 \begin{equation}\label{help2}
 \begin{split}
 \sum_{i=m_1 +1}^{l_2- 1}\sum_{d|(n, i)} \mu(d)
 2^{\lfloor \frac{m_1}{d} \rfloor +\lfloor \frac{m_2}{d} \rfloor - \frac{i}{d}}
  &=
 \sum_{d|n}\sum_{\substack{i=m_1+1\\ d|i}}^{l_2-1}
 \mu(d) 2^{\lfloor \frac{m_1}{d} \rfloor +\lfloor \frac{m_2}{d} \rfloor - \frac{i}{d}} \\
 &=
 \sum_{d|n}\mu(d) 2^{\lfloor \frac{m_1}{d} \rfloor +\lfloor \frac{m_2}{d} \rfloor}
 \sum_{j=\lfloor \frac{m_1}{d} \rfloor +1}^{\lfloor \frac{l_2-1}{d} \rfloor} 2^{-j} \\
 &=
 \sum_{d|n}\mu(d) 2^{\lfloor \frac{m_2}{d} \rfloor}
 \left(1- 2^{-\lfloor \frac{l_2-1}{d}\rfloor +\lfloor \frac{m_1}{d}\rfloor} \right).
 \end{split}
 \end{equation}
 Now combining identities (\ref{help1}, \ref{help2}) yields the result.

 \noindent
 (b)
 Proceeding as in part (a) we find
  \begin{equation}\label{help3}
  \begin{split}
  \Phi_k ([1,m_1]\cup [l_2,m_2],n) &=
  \sum_{d|n} \mu(d)
   \binom{\lfloor \frac{m_2}{d}\rfloor}{k} -
   \sum_{i=m_1 +1}^{l_2-1}\sum_{d|(n, i)} \mu(d)
  \binom{\lfloor \frac{m_1}{d} \rfloor +\lfloor \frac{m_2}{d} \rfloor - \frac{i}{d}}{k-1}.
  \end{split}
  \end{equation}
 Rearranging the last summation on the right of (\ref{help3}) gives
 \begin{equation}\label{help4}
 \begin{split}
 \sum_{i=m_1 +1}^{l_2 -1}\sum_{d|(n, i)}
  \binom{\lfloor \frac{m_1}{d} \rfloor +\lfloor \frac{m_2}{d} \rfloor - \frac{i}{d}}{k-1}
  &=
  \sum_{d|n}\mu(d)
  \sum_{j=\lfloor \frac{m_1}{d} \rfloor +1}^{\lfloor \frac{l_2-1}{d} \rfloor}
  \binom{\lfloor \frac{m_1}{d} \rfloor +\lfloor \frac{m_2}{d} \rfloor- j}{k-1}\\
  &=
  \sum_{d|n}\mu(d)
  \sum_{i=\lfloor \frac{m_1}{d} \rfloor +\lfloor \frac{m_2}{d} \rfloor-\lfloor \frac{l_2-1}{d} \rfloor}^{
  \lfloor \frac{m_2}{d} \rfloor -1} \binom{i}{k-1} \\
  &=
  \sum_{d|n}\mu(d) \left(
  \binom{\lfloor \frac{m_2}{d} \rfloor}{k}-
  \binom{\lfloor \frac{m_1}{d} \rfloor +\lfloor \frac{m_2}{d} \rfloor-\lfloor \frac{l_2-1}{d} \rfloor}{k}
  \right),
  \end{split}
  \end{equation}
  where the last identity follows by formula (\ref{binomial}). Then
  identities (\ref{help3}, \ref{help4}) yield the desired result.
 \end{proof}
  \noindent
 {\bf Definition 2.}
 Let
 \[
 \begin{split}
 \varepsilon(A,B,n) &= \# \{ X \subseteq B:\ X \not=\emptyset,\ X \cap A= \emptyset,\ \text{and\ }
     \gcd(X,n)=1 \}, \\
  \varepsilon_k(A,B,n) &= \# \{ X \subseteq B:\ \# X = k,\ X \cap A= \emptyset,\ \text{and\ }
     \gcd(X,n)=1 \}.
  \end{split}
 \]
 If $B= [1,n]$ we will simply write $\varepsilon(A,n)$ and $\varepsilon_k(A,n)$ rather than
 $\varepsilon(A,[1,n],n)$ and $\varepsilon_k(A,[1,n],n)$ respectively.

  \noindent
  {\bf Theorem 3.} If $l \leq m < n$,
  then
  \[ \text{(a)\ } \varepsilon([l,m],n) = \sum_{d|n} \mu(d) 2^{\lfloor (l-1)/d \rfloor + n/d - \lfloor m/d \rfloor}, \]
  \[ \text{(b)\ } \varepsilon_k([l,m],n) = \sum_{d|n} \mu(d) \binom{\lfloor (l-1)/d \rfloor + n/d - \lfloor m/d \rfloor}{k} . \]
  \begin{proof}
  Immediate from Theorem 2 since
  \[ \varepsilon([l,m],n) = \Phi([1,l-1]\cup [m+1,n],n)\ \text{and\ }
   \varepsilon_k([l,m],n) = \Phi_k([1,l-1]\cup [m+1,n],n). \]
  \end{proof}
 \section*{\normalsize 3. Relatively prime supersets}
 In this section the sets $A$ and $B$ are not necessary nonempty.

 \noindent

 {\bf Definition 3.}
  If $A\subseteq B$ let
 \[
 \begin{split}
 \overline{\Phi}(A,B,n) &= \# \{X\subseteq B:\ X \not= \emptyset,\ A\subseteq X,\ \text{and\ } \gcd(X,n)=1 \}, \\
 \overline{\Phi}_k(A,B,n) &= \# \{X\subseteq B:\ A\subseteq X,\ \# X=k,\ \text{and\ } \gcd(X,n)=1 \}, \\
 \overline{f}(A,B) &= \# \{X\subseteq B:\ X \not= \emptyset,\ A\subseteq X,\ \text{and\ } \gcd(X)=1 \}, \\
 \overline{f}_k (A,B) &= \# \{X\subseteq B:\ \# X = k,\ A\subseteq X,\ \text{and\ } \gcd(X)=1 \}.
 \end{split}
 \]

 \noindent
 The purpose of this section is to give formulas for $\overline{f}(A,[l,m])$,
  $\overline{f}_k(A,[l,m])$, $\overline{\Phi}(A,[l,m],n)$, and $\overline{\Phi}_k(A,[l,m],n)$ for any subset
 $A$ of $[l,m]$. We need a lemma.

 \noindent
 {\bf Lemma 2.}
  If $A \subseteq [1,m]$, then
 \[ \text{(a)\quad }
 \overline{\Phi}(A,[1,m],n) = \sum_{d| (A,n)}\mu(d) 2^{\lfloor m/d \rfloor - \# A},
 \]
 \[ \text{(b)\ }
 \overline{\Phi}_k(A,[1,m],n) = \sum_{d| (A,n)}\mu(d) \binom{ \lfloor m/d \rfloor - \# A}{k- \# A}\
 \text{whenever\ } \# A \leq k \leq m.
 \]

 \begin{proof}
 If $A = \emptyset$, then clearly
 \[
  \overline{\Phi}(A,[1,m],n) = \Phi([1,m],n) \ \text{and\ }
  \overline{\Phi}_k (A,[1,m],n) = \Phi_k ([1,m],n)
 \]
 and the identities in (a) and (b) follow by Theorem 1 for $l=1$.
 Assume now that $A \not= \emptyset$. If $m\leq n$, then
 \[ 2^{m- \# A} = \sum_{d|(A,n)} \overline{\Phi}(\frac{A}{d},[1,\lfloor m/d \rfloor],n/d) \]
 and
 \[ \binom{m- \# A}{k- \# A} = \sum_{d| (A,n)}\mu(d) \overline{\Phi}_k(\frac{A}{d},[1,\lfloor m/d \rfloor],n/d) \]
 which by  M\"obius inversion
  \cite[Theorem 2]{ElBachraoui1} are equivalent to the identities in (a) and in (b) respectively.
  If $m >n$, let $a$ be a positive integer such that $m \leq n^a$. As
 $\gcd(X,n)=1$ if and only if $\gcd(X,n^a)=1$ and $\mu(d) =0$ whenever $d$ has a nontrivial
 square factor we have
 \[
 \begin{split}
  \overline{\Phi}(A,[1,m],n) &= \overline{\Phi}(A,[1,m],n^a) \\
  &= \sum_{d| (A,n^a)}\mu(d) 2^{\lfloor m/d \rfloor - \# A} \\
  &= \sum_{d| (A,n)}\mu(d) 2^{\lfloor m/d \rfloor - \# A}.
  \end{split}
 \]
 The same argument gives the formula for $\overline{\Phi}_k(A,[1,m],n)$.
 \end{proof}

 \noindent
 {\bf Theorem 4.}\label{thm:main3}
 If $A\subseteq [l,m]$, then
 \[ \text{(a)\quad }
 \overline{\Phi}(A,[l,m],n)= \sum_{d| (A,n)}\mu(d) 2^{\lfloor m/d \rfloor - \lfloor (l-1)/d \rfloor -\# A}, \]
 \[
  \text{(b)\quad }
  \overline{\Phi}_k (A,[l,m],n)= \sum_{d| (A,n)}\mu(d) \binom{ \lfloor m/d \rfloor - \lfloor (l-1)/d \rfloor -\# A}{k- \# A}\
  \text{whenever\ } \# A \leq k \leq m-l+1. \]

 \begin{proof}
 If $A= \emptyset$, then clearly
 \[ \overline{\Phi}(A,[l,m],n)= \Phi ([l,m],n) \] and
 \[ \overline{\Phi}_k (A,[l,m],n)= \Phi_k ([l,m],n) \]
 and the identities in (a) and (b) follow by Theorem 1. \\
 Assume now that $A\not= \emptyset$.
  Let
 \[ \Psi (A,l,m,n) = \#
  \{ X\subseteq [l,m]:\ A\cup\{l\} \subseteq X, \text{and\ } \gcd(X,n)=1 \}.
 \]
 Then
 \[
 2^{m-l- \# A}= \sum_{d|(A,l,n)} \Psi (\frac{A}{d}, l/d, \lfloor m/d \rfloor, n/d), \]
 which by M\"obius inversion \cite[Theorem 2]{ElBachraoui1} means that
 \begin{equation} \label{eq:one}
  \Psi (A,l,m,n) = \sum_{d|(A,l,n)} \mu(d) 2^{\lfloor m/d \rfloor -l/d - \# A}.
 \end{equation}
 Then combining identity (\ref{eq:one}) with Lemma 2 gives
 \begin{equation}
 \begin{split}
  \overline{\Phi}(A,[l,m],n) &= \overline{\Phi}([A,[1,m],n) - \sum_{i=1}^{l-1} \Psi(i,m,A,n) \\
  &= \sum_{d|(A,n)} \mu(d) 2^{ \lfloor m/d \rfloor - \# A} -
    \sum_{i=1}^{l-1} \sum_{d|(A,i,n)} \mu(d) 2^{\lfloor m/d \rfloor -i/d - \# A} \\
  &= \sum_{d|(A,n)} \mu(d) 2^{ \lfloor m/d \rfloor - \# A} -
    \sum_{d|(A,n)} \mu(d) 2^{ \lfloor m/d \rfloor - \# A} \sum_{j=1}^{\lfloor (l-1)/d \rfloor} 2^{-j} \\
  &= \sum_{d|(A,n)} \mu(d) 2^{ \lfloor m/d \rfloor - \# A} -
    \sum_{d|(A,n)} \mu(d) 2^{ \lfloor m/d \rfloor - \# A}(1 - 2^{- \lfloor (l-1)/d \rfloor}) \\
  &= \sum_{d| (A,n)}\mu(d) 2^{\lfloor m/d \rfloor - \lfloor (l-1)/d \rfloor -\# A}.
 \end{split}
 \end{equation}
 This completes the proof of (a). Part (b) follows similarly.
 \end{proof}

 \noindent
 As to $\overline{f}(A,[l,m])$ and $\overline{f}_k (A,[l,m])$ we similarly have:

 \noindent
 {\bf Theorem 5.}\label{thm:main4}
 If $A \subseteq [l,m]$, then

 \[ \text{(a)\quad }
 \overline{f}(A,[l,m])= \sum_{d| \gcd(A)}\mu(d) 2^{\lfloor \frac{m}{d} \rfloor - \lfloor \frac{l-1}{d} \rfloor -\# A}, \]
 \[ \text{(b)\quad }
  \overline{f}_k (A,[l,m])= \sum_{d| \gcd(A)}\mu(d) \binom{ \lfloor \frac{m}{d}\rfloor - \lfloor \frac{l-1}{d} \rfloor -\# A}{k- \# A},\
  \text{whenever\ } \# A \leq k\leq m-l+1. \]

 We close this section by formulas for relatively prime sets which have a nonempty intersection
 with $A$.

 \noindent
 {\bf Definition 4.}
 Let
 \[
 \begin{split}
  \overline{\varepsilon}(A,B,n) &= \# \{ X \subseteq B:\ X \cap A\not= \emptyset\ \text{and\ }
     \gcd(X,n)=1 \}, \\
  \overline{\varepsilon}_k(A,B,n) &= \# \{ X \subseteq B:\ \# X =k,\ X \cap A\not= \emptyset,\ \text{and\ }
     \gcd(X,n)=1 \}, \\
  \overline{\varepsilon}(A,B)  &= \# \{ X \subseteq B:\ X \cap A\not= \emptyset\ \text{and\ }
     \gcd(X)=1 \}, \\
  \overline{\varepsilon}_k(A,B)  &= \# \{ X \subseteq B:\ \# X =k,\ X \cap A\not= \emptyset,\ \text{and\ }
     \gcd(X)=1 \}.
  \end{split}
  \]

 \noindent
 {\bf Theorem 6.}
 We have
 \[
  \text{(a)\quad } \overline{\varepsilon}(A,[l,m],n)= \sum_{\emptyset\not= X \subseteq A}
  \sum_{d|(X,n)} \mu(d) 2^{\lfloor \frac{m}{d} \rfloor - \lfloor \frac{l-1}{d} \rfloor -\# X},
 \]
 \[
 \text{(b)\quad } \overline{\varepsilon}_k(A,[l,m],n)= \sum_{\substack{\emptyset \not= X\subseteq A \\ \# X \leq k}}
  \sum_{d|(X,n)} \mu(d) \binom{ \lfloor \frac{m}{d}\rfloor - \lfloor \frac{l-1}{d} \rfloor -\# X}{k- \# X},
  \]
  \[
  \text{(c)\ } \overline{\varepsilon}(A,B)  = \sum_{\emptyset\not= X \subseteq A}
  \sum_{d|\gcd(X)} \mu(d) 2^{\lfloor \frac{m}{d} \rfloor - \lfloor \frac{l-1}{d} \rfloor -\# X},
 \]
 \[
 \text{(d)\ } \overline{\varepsilon}_k(A,B)  = \sum_{\substack{\emptyset \not= X\subseteq A \\ \# X \leq k}}
  \sum_{d|\gcd(X)} \mu(d) \binom{ \lfloor \frac{m}{d}\rfloor - \lfloor \frac{l-1}{d} \rfloor -\# X}{k- \# X}.
  \]

  \begin{proof}
  These formulas Follow by Theorems 4, 5 and the facts that
  \[ \overline{\varepsilon}(A,[l,m],n)= \sum_{\emptyset\not= X \subseteq A} \overline{\Phi}(X,[l,m],n), \]
  \[ \overline{\varepsilon}_k(A,[l,m],n)= \sum_{\substack{\emptyset \not= X\subseteq A \\ \# X \leq k}}
  \overline{\Phi}_k (X,[l,m],n),
  \]
  \[ \overline{\varepsilon}(A,[l,m])  = \sum_{\emptyset\not= X \subseteq A}\overline{f}(X,[l,m]), \]
  \[ \overline{\varepsilon}_k(A,[l,m])  = \sum_{\substack{\emptyset\not= X \subseteq A \\ \# X \leq k}}
  \overline{f}_k(X,[l,m]). \]

  \end{proof}


\end{document}